\documentclass[a4paper,12pt]{article}
\usepackage[fleqn]{amsmath}
\usepackage{blkarray}
\usepackage{multirow}
\usepackage{mathtools}
 \usepackage{graphicx}
\usepackage[width=17cm,height=25cm]{geometry}
 \usepackage{epsfig}
\usepackage{amssymb}
\usepackage{amsmath}

\newcommand{\be}{\begin{equation}}
\newcommand{\ee}{\end{equation}}
\newcommand{\bea}{\begin{eqnarray}}
\newcommand{\eea}{\end{eqnarray}}
\newcommand{\bbea}{\begin{eqnarray*}}
\newcommand{\eeea}{\end{eqnarray*}}

\newtheorem{theorem}{Theorem}[section]
\newtheorem{lemma}{Lemma}[section]
\newtheorem{ass}{Assumption}[section]
\newtheorem{example}{Example}[section]
\newtheorem{definition}{Definition}[section]
\newtheorem{proposition}{Proposition}[section]
\newtheorem{corollary}{Corollary}[section]

\newcommand{\gy}{\psi}

\DeclareMathOperator*{\slim}{s-lim}

\begin{document}

\title{Spectral and scattering theory for perturbed block Toeplitz operators}

\author{Petru Cojuhari\footnote{Faculty of Applied Mathematics, AGH University of Science and Technology, al. Mickiewicza 30, 
30-059  Krakow,  Poland.  cojuhari@agh.edu.pl }  
and Jaouad Sahbani\footnote{Institut de Math\'ematiques de Jussieu-Paris Rive Gauche-UMR7586,
 Universit{\'e} Paris Diderot, B\^atiment Sophie Germain--case 7012, 5 rue Thomas Mann, 75205 Paris Cedex 13, FRANCE.
 jaouad.sahbani@imj.prg.fr}}
\date{\today}

\maketitle
\begin{abstract}
We analyse  spectral properties of a class of compact perturbations of  block Toeplitz operators associated with analytic symbols.
In particular,  a limiting absorption principle  and  the absence of singular continuous spectrum are shown.
The existence and the completeness of wave operators are also obtained.
Our study is based on the construction of a conjugate operator in Mourre sense for the corresponding Laurent operators. 
\end{abstract}




\section{Introduction and main results}
Let $\mathbb{Z}_+$ be the set of nonnegative integers and
  ${\mathcal H}_+=l^2(\mathbb{Z}_+, \mathbb{C}^N)$, for some integer  $N\geq1$, be the Hilbert
space of square summable vector-valued sequences $ (\gy_n)_{n\geq 0} $ endowed
with the scalar product
$$
\langle\phi,\psi\rangle=\sum_{n=0}^{\infty}\langle\phi_n,\psi_n{\rangle}.
 $$
Let  $(A_j)_{j\in\mathbb{Z}}$ be a family of square $N\times N$ matrices such that
 \begin{equation}\label{hypothese0}
\sum_{j=-\infty}^{+\infty}\|A_j\|<\infty\quad\mbox{ and } \quad   A_{-j}=A_j^*, \quad  j\in\mathbb{Z}.
\end{equation}
Here-above we denoted by $A^*$ the adjoint matrix of a given matrix $A$.
 Consider the Toeplitz operator
$T_0$  defined in ${\mathcal H}_+$  by the expression
 \begin{equation}\label{defH0}
(T_0\psi)_n= \sum_{m=0}^{+\infty}A_{n-m}\psi_m, ~~\mbox{\quad} n=0,1,2,....
 \end{equation}
According to (\ref{hypothese0}),  $T_0$ is a bounded self-adjoint operator  in ${\mathcal H}_+$. 
Our goal here is to study the spectral properties of $T_0$ and its compact perturbations.
The crucial step of our analysis is the construction of a conjugate operator  in the sense of  the Mourre estimate \cite{M} (see also Section \ref{rappel})
 for a large class of compact perturbations of  the corresponding Laurent operator.  

We will need the following standard notations. 
For a self-adjoint operator $S$ we denote by $E_S(\cdot)$ its spectral measure, $\sigma(S)$ its spectrum,
$\sigma_{ess}(S)$ its essential spectrum,
 $\sigma_p(S)$ the set of its eigenvalues,  $\sigma_{sc}(S)$ its singular continuous spectrum, and
  $\sigma_{ac}(S)$ its absolutely continuous spectrum. The imaginary part of the complex number $z$ will be denoted by $\Im(z)$.

Consider on the unit circle
$\mathbb{U}=\{e^{ip}\in\mathbb{C}~/~p\in\mathbb{R}\}$  the Hermitian matrix-valued function  defined by
\begin{equation}\label{reduced0}
h(e^{ip})=\sum_{m=-\infty}^{\infty}A_me^{imp}, \quad p\in\mathbb{R}.
\end{equation}
It is usually called in the literature \cite{GK} (see also \cite{GF,BoSi}) the symbol of $T_0$ and
plays a crucial role in its spectral analysis.
According to (\ref{hypothese0}), $h$  is continuous, but
here we need to assume the following.
\begin{ass}\label{analytique}
 The Hermitian matrix-valued function $h$ is holomorphic on unit circle, that is, it has an
 holomorphic extension to  an open annulus about the unit circle $\mathbb{U}$ of the form 
 $\{z\in\mathbb{C}~/1/r<|z|<r\}$ for some $r>1$.
 \end{ass}
In the sequel, we  identify $h$ with  the $2\pi-$periodic function  
$
p\mapsto h(e^{ip})
$
 that we denote by the same symbol $h$.
A simplest example is obtained when the sequence $A_j=0$ as soon as $|j|>M$, for some integer $M$. In this case $h$ is a trigonometric polynomial and $T_0$
is a block banded matrix. The case  $M=1$ leads to block Jacobi matrices.

Under  Assumption \ref{analytique}, see \cite{Ka} and also \cite{Ba},  there exist N real-valued $2\pi-$periodic functions
$\{\lambda_{j}(p);\,1\leq j\leq N\}$ representing the repeated eigenvalues  of
$h(p)$ for all $p\in[-\pi,\pi)$, 
and N vector-valued $2\pi-$periodic functions $\{ W_{j}(p),\,1\leq j\leq N\}$ representing a corresponding orthonormal basis
of eigenvectors. These functions are  smooth (they are even analytic) on $(-\pi,\pi)$ and they are piecewise smooth on $\mathbb{R}$. 
In particular,  for any $j=1,\cdots,N$,  the set
$$
\kappa_1(\lambda_j)=\{\lambda_j(p)/\lambda_j \mbox{ is not diffrentiable at $p$}\}
$$
contains at most two values.
Therefore,  for any $j=1,\cdots,N$,  the set of critical values of $\lambda_j$  defined by 
$$
\kappa(\lambda_j)=\kappa_1(\lambda_j)\cup\{ \lambda_j(p)~/~p\in(-\pi,\pi)\mbox{ and }  \lambda_j'(p)=0\}
$$
 is clearly  finite. 
Finally, the critical set of $h$ defined by
$$
\kappa(h)=\cup_{i=1}^{N}\kappa(\lambda_j)
$$
is finite, too.  It is known that  \cite{GK}  the essential spectrum $\sigma_{ess}(T_0)$ of $T_0$ consists of the union of  $N$ compact intervals, called spectral bands of $T_0$.
More precisely,  let $\alpha_j=\min_{p\in[-\pi,\pi]}\lambda_j(p)$ and 
$\beta_j=\max_{p\in[-\pi,\pi]}\lambda_j(p)$ for each $j=1,\cdots, N$. Then
$$
\sigma_{ess}(T_0)=\cup_{j=1}^N\Sigma_j\quad \mbox{ with } \quad \Sigma_j=[\alpha_j,\beta_j], \mbox{  } ~ j=1,\cdots,N, 
$$
Some of these spectral bands  $\Sigma_j$ may degenerate into single point, i.e., $\alpha_{j}=\beta_j$, 
in which case, this value is an infinitely degenerate eigenvalue of $T_0$ and belongs to $\kappa(h)$. 
For example, if
$
h(p)=\left(
\begin{array}{cc}
  0&e^{-ip}    \\
  e^{ip} &  0   
\end{array}
\right),
$
then clearly
$
\sigma_{ess}(T_0)=\kappa(h)=\{-1,1\}.
$
 Usually in the literature one avoids this degeneracy by assuming that for any scalar $\lambda$ the function 
 $\det(h(p)-\lambda)$ does not vanish identically  with respect to $p$.
We will not do that since our approach allows us to study what happens outside the set of critical values $\kappa(h)$.
There is also cases where some spectral bands of $T_0$ are non degenerate and some are degenerate. For example, 
 if
$
h(p)=\left(
\begin{array}{cc}
  \cos p&0    \\
  0 &  0   
\end{array}
\right),
$
then 
$$
\Sigma_1=\sigma_{ess}(T_0)=[-1,1], ~~\Sigma_2=\kappa(h)=\{0\},
$$ 
and 0 is an infinitely degenerate eigenvalue of $T_0$.

Recall that in the scalar case, i.e. $N=1$,  $T_0$ has no discrete spectrum, thanks to
 Hartman-Wintner's theorem \cite{HW}, it follows, also from \cite{K}, and we refer to \cite{BoSi,D,GF}  for more information on the subject. 
In contrast, if $N>1$ then the discrete spectrum $\sigma_d(T_0)$ of $T_0$ is non empty in general. 
For example, if
$$
h(p)=\left(
\begin{array}{cc}
  0&ae^{-ip}+b    \\
  ae^{ip}+b &  0   
\end{array}
\right) ~~ \quad\mbox{with $a,b>0$.}
$$
then one may show that
the essential spectrum of the associated Toeplitz operator  $T_0$  is 
$$
\sigma_{ess}(T_0)=\left[-(a+b),-|b-a|\right]\cup\left[ |b-a|,a+b \right].
$$
Here the spectral bands are separated by a non trivial gap given by $(-|b-a|,|b-a|)$  if and only if $a\not= b$.
Moreover, if $a>b$ then 0 is an eigenvalue of $T_0$ and  a corresponding eigenvector $\psi$ is given by $\psi_n=((-b/a)^n,0),n\geq0$.
In contrast,  if $a=b$ then the eigenvalues of $h(p)$ are given on $[-\pi,\pi)$ by
$$
\lambda_\pm(p)=\pm2a\cos(p/2)\quad\mbox{and } \quad \sigma_{ess}(T_0)=[-2a, 2a].
$$
Observe that these eigenvalues are only piecewise smooth as $2\pi$-periodic functions on $\mathbb{R}$. 
This is because  $p=-\pi$ is an exceptional point where the eigenvalues of $h(p)$ meet each other.
Mention also that,  under Assumption \ref{analytique}, the discrete spectrum  of $T_0$ in each spectral gap is only finite. It was proved in \cite{F} (based on an abstract result from \cite{CF}) and also in \cite{C1}.
Moreover, in \cite{CR,C1} this fact is discussed even for the point spectrum.
Finally, the absence of singular continuous spectrum of $T_0$ is proved in \cite{R} in the case where the symbol $h$ is a trigonometric polynomial.
Here we will prove this assertion under Assumption \ref{analytique}  by a different method based on the positive commutator approach.
 Moreover, we obtain sharp information on the
 asymptotic behaviour of the resolvent of  $T_0$  near the real axis, see our next Theorem \ref{thm2}.
 It is worth to add that our approach can be extended to a large class of smooth symbols $h$ in a straightforward manner.
For a similar analysis in the scalar case $N=1$ one may see \cite{ABC,BS1}. 

\begin{theorem}\label{thm1}
\begin{enumerate}
 \item  Outside $\kappa(h)$ the  eigenvalues of $T_0$  are all finitely degenerate 
and their possible accumulation points are included in $\kappa(h)$. In particular, the set 
$$
\tau(T_0):=\sigma_p(T_0)\cup\kappa(h)
$$
 is closed and countable real subset.
\item The singular continuous spectrum of $T_0$ is empty.
\end{enumerate}
\end{theorem}
To show the absence of the singular continuous spectrum  of $T_0$ we actually establish a limiting absorption principle for $T_0$.
In order to state this we need the following notations.
Let $\textbf{N}$ be  the multiplication operator defined by 
$(\textbf{N}\psi)_j=j\psi_j$ for all $\psi\in\mathcal{H}_+$
with its natural domain $D(\textbf{N})=\{\psi\in\mathcal{H}_+~/~\textbf{N}\psi\in\mathcal{H}_+\}$.
 We denote  by  $B(\mathcal{E},\mathcal{F})$ the space of bounded operators between the normed vector spaces $\mathcal{E}$ and $\mathcal{F}$.
If $\mathcal{E}=\mathcal{F}$ then we put $B(\mathcal{E})=B(\mathcal{E},\mathcal{E})$. 
We also need the interpolation space $\mathcal{K}_+:=(D(\textbf{N}),\mathcal{H}_+)_{1/2,1}$ which can be described, according to  Theorem 3.6.2 of  \cite{ABG}, by the norm, 
\begin{equation}\label{besov}
\|\psi\|_{1/2,1}=\sum_{j=0}^\infty2^{j/2}\|\theta(2^{-j}|\textbf{N}|)\psi\|,
\end{equation}
where $\theta\in C^\infty_0(\mathbb{R})$ with  $\theta(x)>0$ if $1<x<2$ and $\theta(x)=0$ otherwise. Let us denote by $\mathcal{K}_+^*$ the dual space of $\mathcal{K}_+$.
Hence by Riesz identification of $\mathcal{H}_+$ with its dual we have the continuous inclusions
 $\mathcal{K}_+\subset\mathcal{H}_+\subset\mathcal{K}_+^*$. Let us define $\mathcal{K}_+^{*\circ}$ as the closure of $\mathcal{H}_+$ in  $\mathcal{K}_+^*$. 
\begin{theorem}\label{thm2}
\begin{enumerate}
\item The holomorphic maps defined on $\mathbb{C}_\pm=\{z\in\mathbb{C}~/~\pm\Im(z)>0\}$ by $z\mapsto(T_0-z)^{-1}\in B(\mathcal{K}_+,\mathcal{K}_+^*)$ extend weak*-continuously to 
 $\mathbb{C}_\pm\cup\mathbb{R}\setminus\tau(T_0)$. 
 \item  For any bounded operator $B:\mathcal{K}_+^{*\circ}\to\mathcal{F}$, $\mathcal{F}$ a given Hilbert space, 
then $B$ is  locally $T_0-$smooth on $\mathbb{R}\setminus\tau(T_0)$, i.e., for any compact $K\subset\mathbb{R}\setminus\tau(T_0)$ there is $C>0$ such that
$$
\int_{-\infty}^{\infty}\|Be^{-iT_0t}f\|_{\mathcal{F}}^2dt\leq C\|f\|^2~,\quad f\in E_{T_0}(K)\mathcal{H}_+.
$$
\end{enumerate}
 \end{theorem}
In particular, the  limits 
$
 ({T_0}-x\mp i0)^{-1}:=\lim_{\mu\rightarrow0+}({T_0}-x\mp i\mu)^{-1}
$
exist locally uniformly on  $\mathbb{R}\setminus\tau(T_0)$ for the weak* topology of 
$B({\mathcal{K}_+,\mathcal{K}_+^*})$. Actually we will prove also that,  for any $s>1/2$, the maps 
$
 x\mapsto\langle {\textbf{N}}\rangle^{-s}({T}_0-x\mp i0)^{-1}\langle {\textbf{N}}\rangle^{-s}\in B({\mathcal   H}_+)
 $
are locally of class $\Lambda^{s-\frac{1}{2}}$ on $\mathbb{R}\setminus\tau(T_0)$ (see next section for the definition of $\Lambda^\alpha$). 
Here for any real number $x$ we put $\langle {x}\rangle=\sqrt{1+x^2}$ .
%
%
%
%
%
\begin{corollary}\label{corollaire1}
Assume that $A_j=0$ if $|j|>M$ for some integer $M\geq1$, and denote $T_0$ the associated Toeplitz operator.
Then the assertions of Theorems \ref{thm1} and \ref{thm2} hold true. Moreover, the discrete spectrum of 
 $T_0$ is at most  finite.
\end{corollary}
Now,  it is natural to study the  preservation of the spectral structure of $T_0$ under compact perturbations.
More precisely, let $v$ be a symmetric compact operator  in $\mathcal{H}_+$. Then $T=T_0+v$ is a self-adjoint operator 
in $\mathcal{H}_+$ and, thanks to Weyl criterion, 
$
\sigma_{ess}(T)=\sigma_{ess}(T_0).
$
Moreover,  since $v$ is compact,  it is easy to show that 
$$
\lim_{r\to\infty}\|\theta(\langle \mathbf{N}\rangle/r) v\|=0,
$$
where $\theta\in C^\infty_0(\mathbb{R})$ can be taken as in (\ref{besov}). In fact, conditions on how fast this convergence takes place will be given 
to control the singular continuous spectrum.
\begin{theorem}\label{toeplitz-perturbed01}
Let $v$ be a compact symmetric operator in $\mathcal{H}_+$ such that
\begin{equation}\label{admissible1}
\int_1^{\infty}\|\theta(\langle \mathbf{N}\rangle/r) v\|dr<\infty.
\end{equation}
Then
 \begin{enumerate}
 \item outside  $\kappa(h)$ the eigenvalues of $T=T_0+v$ are all finitely degenerate and 
may accumulate only to points in $\kappa(h)$. Hence  $\tau(T):=\kappa(h)\cup\sigma_p(T)$
 is closed and countable;
 \item the holomorphic maps defined on $\mathbb{C}_\pm=\{z\in\mathbb{C}~/~\pm\Im(z)>0\}$ by $z\mapsto(T-z)^{-1}\in B(\mathcal{K}_+,\mathcal{K}_+^*)$ extend weak*-continuously to 
 $\mathbb{C}_\pm\cup\mathbb{R}\setminus\tau(T)$. In particular, $T$ has no singular continuous spectrum.
 \item
For any bounded operator $B:\mathcal{K}_+^{*\circ}\to\mathcal{F}$, $\mathcal{F}$ a given Hilbert space, 
then $B$ is  locally $T-$smooth on $\mathbb{R}\setminus\tau(T)$, i.e., for any compact $K\subset\mathbb{R}\setminus\tau(T)$ there is $C>0$ such that
$$
\int_{-\infty}^{\infty}\|Be^{-iTt}f\|_{\mathcal{F}}^2dt\leq C\|f\|^2~,\quad f\in E_{T}(K)\mathcal{H}_+.
$$
\end{enumerate}
 \end{theorem}
Here is a notable consequence of Theorems \ref{thm1} and  \ref{toeplitz-perturbed01} on the scattering theory of $T$ and $T_0$.
\begin{corollary}\label{waveoperators}
In addition of assumptions of Theorem \ref{toeplitz-perturbed01} assume that  $v$ extends to a bounded operator from $\mathcal{K}_+^{*\circ}$ to 
 $\mathcal{K}_+$. 
 Let $P_{ac}(T_0)$ be the spectral projector of $T_0$ on its purely absolutely continuous component. 
 Then the wave operators
 $$
 \Omega^\pm(T,T_0)=\slim_{t\rightarrow\pm\infty}e^{iTt}e^{-iT_0t}P_{ac}(T_0)
 $$
exist and they are complete, i.e. their range equal the absolutely continuous subspace  $\mathcal{H}_{ac}(T)$ of $T$.
\end{corollary}
 Notice that obviously our condition  (\ref{admissible1}) follows if for some $s>1/2$, 
 \begin{equation}\label{admissible2}
\sup_{r\geq 1}\|r^{s+\frac{1}{2}}\theta(\langle \mathbf{N}\rangle/r) v\|<\infty.
\end{equation}
But under this assumption we also have the following.
\begin{theorem}\label{toeplitz-perturbed02}
Let $v$ be a compact symmetric operator in $\mathcal{H}_+$ such that
 (\ref{admissible2}) holds. Then Theorem \ref{toeplitz-perturbed01}. Moreover, the maps 
$
 x\mapsto\langle {\textbf{N}}\rangle^{-s}({T}-x\mp i0)^{-1}\langle {\textbf{N}}\rangle^{-s}\in B({\mathcal   H}_+)
 $
are locally of class $\Lambda^{s-\frac{1}{2}}$ on $\mathbb{R}\setminus\tau(T)$.
 \end{theorem}

Let $B=(b_{ij})_{i,j\in\mathbb{Z}}$ be a bounded symmetric operator in $\mathcal{H}_+$. Put 
\begin{eqnarray}
n_B(r)&:=&(\sup_{r\leq i\leq 2r }\sum_{j\geq0}\|b_{ij}\|)\cdot(\sup_{ j\geq0}\sum_{r\leq i\leq 2r}\|b_{ij}\|)\\
p_B(r)&:=&\sum_{r\leq i\leq 2r} \sum_{ j\geq0}\|b_{ij}\|^2.
\end{eqnarray}
\begin{corollary}
If $v=v_1+v_2$ such that 
$$
\int_1^{\infty}\sqrt{n_{v_1}(r)}+\sqrt{p_{v_{2}}(r)}dr<\infty,
$$
then (\ref{admissible1}) holds. 
\end{corollary}
In particular, Theorems  \ref{waveoperators} and  \ref{toeplitz-perturbed01} follow. Similarly, condition (\ref{admissible2}) and so Theorem \ref{toeplitz-perturbed02}  hold if,  for some $s>1/2$,
$$
\sup_{r\geq1}r^{s+\frac{1}{2}}(\sqrt{n_{v_1}(r)}+\sqrt{p_{v_{2}}(r)})<\infty,
$$ 

\begin{example}
Let $A_n$ and $B_n$ be two  sequences of $N$ by $N$ matrices such that $B_n=B_n^*$ for all $n\in\mathbb{Z}_+$. 
 Consider  the block Jacobi operator $J=J(\{A_n\},\{B_n\})$ acting in $\mathcal{H}^+$  by
\begin{equation}\label{Jacobi1}
(J\psi)_n=A_{n-1}^*\psi_{n-1}+B_n\psi_n+A_{n}\psi_{n+1} ~~\mbox{
for all } n\geq0
\end{equation}
initialized with $\psi_{-1}=0$.
Assume that there exist two matrices $A$ and $B$ such that
\begin{equation}\label{Jacobi2}
\lim_{n\rightarrow\infty}(\|A_n-A\|+\|B_n-B\|)=0.
\end{equation}
So that $J=T_0+v$, where $T_0$ is the Toeplitz operator associated with the symbol
\begin{equation}\label{symbol1}
h(p)=e^{-ip}A+e^{ip}A^*+B.
\end{equation}
Then clearly our results applies directly if the convergence in (\ref{Jacobi2}) is sufficiently fast. For example,
outside $\kappa(h)$  the eigenvalues of $J$ are all finitely degenerate and their possible accumulation points are contained  in $\kappa(h)$ as soon  as
\begin{equation}\label{HP1}
\lim_{n\rightarrow\infty}n(\|A_n-A\|+\|B_n-B\|)=0.
\end{equation}
Similarly, conditions of Theorem \ref{toeplitz-perturbed01}  and  Corollary \ref{waveoperators} are satisfied 
 if 
 $$
 \int_1^{\infty}\sup_{r\leq n\leq 2r }(\|A_n-A\|+\|B_n-B\|)dr<\infty.
 $$
\end{example}

The paper is organized as follows.   Section \ref{rappel} contains a brief review on what we need from Mourre's theory.  
In Section \ref{LO} we study in details the Laurent operator $H_0$ associated with the symbol $h$. In particular,
a Mourre estimate is established locally for $H_0$ outside the critival set $\kappa(h)$. 
 In section \ref{perturbŽ} we  extend our Mourre estimate for a large class of compact perturbations of $H_0$. 
 In section \ref{toelplitz-analytique} Theorems \ref{thm1} and \ref{thm2} are proven. Finally, we show
 Theorems \ref{toeplitz-perturbed01} and \ref{toeplitz-perturbed02}
 in Section  \ref{Toeplitz-pertubation} .

\protect\setcounter{equation}{0}
%
%
%
\section{ The conjugate operator theory}\label{rappel}
The following brief review on the conjugate operator theory is based on
\cite{ABG,BG,BGS2,S1}. 
Let $(E,\|\cdot\|)$ be a Banach space and $f:{\Bbb R}\rightarrow E$ be a bounded
continuous function.  For an integer $m\geq1$ let $w_m$ be the modulus of
continuity of order $m$ of $f$ defined on $(0,1)$ by
$$
w_m(f,\varepsilon)=\sup_{x\in\mathbb{R}}
\biggl\|\sum_{j=0}^m (-1)^j
\left(
\begin{array}{cc}
 m  \\ j   
\end{array}
\right)
f(x+j\varepsilon)\biggr\|.
$$
We say that $f\in\Lambda^{\alpha}$, $\alpha>0$,
 if  there is an integer $m>\alpha$ such
that 
$$
\sup_{0<\varepsilon<1}\varepsilon^{-\alpha} w_m(f,\varepsilon)<\infty.
$$
Notice that if $\alpha\in(0,1)$, then $\Lambda^\alpha$ is nothing but the space of H\H{o}lder continuous functions of order $\alpha$.
In contrast, if $\alpha=1$ then $\Lambda^1$ consists of smooth functions in Zygmund's sense (they are not Lipschitz in general, see \cite{Z}). Finally, if $\alpha>1$ and 
$n_\alpha$ is the greatest integer strictly less than $\alpha$,  then $f$ belongs to $\Lambda^\alpha$ if and only if $f$ is $n_\alpha-$times continuously differentiable with bounded derivatives and its derivative
 $f^{(n_\alpha)}$ of order $n_\alpha$  is of class $\Lambda^{\alpha-n_\alpha}$. 
 For example, $f$ belongs to $\Lambda^2$  means that $f$ is continuously differentiable with a bounded derivative and $f'$ is of class  $\Lambda^1$ ($f'$ is not Lipschitz in general).

Let $\mathbb{A}$ be a  self-adjoint operator in a  separable complex Hilbert space  $\mathcal{H}$ and $S\in B({\mathcal H})$.
\begin{definition} \label{defreg}
(i) Let $k\geq 1$ be an integer and $\sigma>0$.
 We say that  $S$ is of class $C^{k}(\mathbb{A})$, respectively of class $C_u^{k}(\mathbb{A})$, $\mathcal{C}^{\sigma}(\mathbb{A})$, if the map
$$
t\longmapsto S(t)= e^{-i\mathbb{A}t}Se^{i\mathbb{A}t}\in B({\mathcal H})
$$
is strongly of class $C^k$, respectively of class $C^k$ in norm, of class $\Lambda^\sigma$ on ${{\Bbb R}}$.

(ii)  We say that $S$ is of class $\mathcal{C}^{1,1}(\mathbb{A})$, if
$$
\int_0^1\|e^{-i\mathbb{A}\varepsilon}Se^{i\mathbb{A}\varepsilon}-2S+e^{i\mathbb{A}\varepsilon}Se^{-i\mathbb{A}\varepsilon}S\|\frac{d\varepsilon}{\varepsilon^2}<\infty.
$$
\end{definition}
\textbf{Remarks} 
(i)  We have the following inclusions 
$$
\mathcal{C}^s(\mathbb{A})\subset \mathcal{C}^{1,1}(\mathbb{A})\subset C_u^1(\mathbb{A})\subset C^1(\mathbb{A}),\quad s>1.
$$

(ii) 
One may show  that $S$ is  of class $C^1(\mathbb{A})$ if and only if the sequilinear form defined on 
$D(\mathbb{A})$ by 
$[S,\mathbb{A}]=S\mathbb{A}-\mathbb{A}S$  has a continuous extension to $\mathcal H$, which we identify with
the associated bounded operator in $\mathcal H$ (from the Riesz Lemma) that
we denote  by the same symbol.  Moreover,
$
[S,i\mathbb{A}]=\frac{d}{dt}|_{t=0}S(t).
$

(iii) Recall that  $\langle\mathbb{A}\rangle=(1+\mathbb{A}^2)^{1/2}$. To prove that $S$  is of class 
$\mathcal{C}^{\sigma}(\mathbb{A})$ it is enough to show that
  $\langle\mathbb{A}\rangle^\sigma S$ is bounded in $\mathcal H$,
see for example the appendix of \cite{BS1}.  In particular,  in the case where $\sigma>1$, the operator $S$ is of class ${C}^{\sigma}(\mathbb{A})$ if one of the following conditions holds true\\
(1) $\langle\mathbb{A}\rangle^\sigma S$ is bounded in $\mathcal H$ or\\
(2) $S$ is of class  ${C}^{1}(\mathbb{A})$ and $\langle\mathbb{A}\rangle^{\sigma-1} [S,i\mathbb{A}]$ is bounded in $\mathcal H$.

In the sequel of this section, let   $H$ be a bounded self-adjoint operator which is at least of class $C^1(\mathbb{A})$. 
Then $[H,i\mathbb{A}]$ defines  a   bounded operator
 in $\mathcal{H}$ that we still denote by the
 same symbol $[H,i\mathbb{A}]$.

 Let  $\Delta$ be a compact interval such that,
for some constant $a>0$ and  a compact operator $K$ in $\mathcal H$, we have
\bea\label{mourre estimate}
E_H(\Delta )[H,i\mathbb{A}]E_H(\Delta )\geq aE_H(\Delta )+K.
\eea
The inequality (\ref{mourre estimate}) is called {\it the Mourre estimate} and we say that $\mathbb{A}$ is conjugate to $S$ on $\Delta$. If   (\ref{mourre estimate}) holds with $K=0$ then we say that $\mathbb{A}$ is strictly conjugate to $S$ on $\Delta$.
\begin{theorem} \label{virial} 
The interval $\Delta$ contains at most a finite number of eigenvalues of $H$ counted with multiplicities. 
Moreover,  if   $K=0$, then the operator $H$ has no eigenvalues in $\Delta$.
\end{theorem}

Let $\mathcal{K}_\mathbb{A}:=\mathcal{H}_{1/2,1}$ be the Besov space associated to 
$\mathbb{A}$ defined by the norm
$$
\|\psi\|_{1/2,1}=\|\tilde\theta(\mathbb{A}) \psi\|+\sum_{j=0}^\infty2^{j/2}\|\theta(2^{-j}|\mathbb{A}|)\psi\|,
$$
where $\theta\in C^\infty_0(\mathbb{R})$ with  $\theta(x)>0$ if $1<x<2$ and $\theta(x)=0$ otherwise, and $\tilde\theta\in C^\infty_0(\mathbb{R})$ with  $\tilde\theta(x)>0$ if $|x|<2$ and $\tilde\theta(x)=0$ otherwise.
One has, see \cite{BG,S1}, the following.
\begin{theorem} \label{lap}
 Assume that $H$ is of class $\mathcal{C}^{1,1}(\mathbb{A})$. Then $H$ has no singular continuous spectrum in $\Delta$. 
 Moreover, the  limits
$
({H}-x\mp i0)^{-1} :=
\lim_{\mu\rightarrow0+}({H}-x\mp i\mu)^{-1}
$
exist locally uniformly on  $\Delta\setminus\sigma_p(H)$ for the weak* topology of $B(\mathcal{K}_\mathbb{A},\mathcal{K}_\mathbb{A}^*)$.
\end{theorem}
Next theorem describes continuity properties of these boundary values of the resolvent of $H$ 
as functions of $x$, as well as some of their propagation consequences. \begin{theorem} \label{lap1}
 Assume that $H$ is of class $\mathcal{C}^{s+\frac{1}{2}}(\mathbb{A})$
for some $s>1/2$. Then the  maps
$
 x\mapsto\langle\mathbb{A}\rangle^{-s}({H}-x\mp i0)^{-1}\langle\mathbb{A}\rangle^{-s}
$
are locally of class
$\Lambda^{s-\frac{1}{2}}$ on $\Delta\setminus\sigma_p(H)$.
Moreover, for every $\varphi\in C^\infty_0(\Delta\setminus\sigma_p(H))$ we have
$$
||\langle\mathbb{A}\rangle^{-s}e^{-iHt}\varphi(H)\langle\mathbb{A}\rangle^{-s}||\leq C<t>^{-(s-\frac{1}{2})}.
$$
\end{theorem}
Since $\|e^{-iHt}\varphi(H)\|=\|\varphi(H)\|$, then by using  the complex interpolation one  obtains much more than this, see \cite{BG,BGS2,S1}.
For example, if $H$ is of class $C^\infty(\mathbb{A})$, then, for any  $\varphi\in C^\infty_0(\Delta\setminus\sigma_p(H))$ and  $\sigma,\varepsilon>0$, we have
$$
||\langle\mathbb{A}\rangle^{-\sigma}e^{-i{H}t}\varphi(H)\langle\mathbb{A}\rangle^{-\sigma}||\leq C(1+|t|)^{-\sigma+\varepsilon},\quad \mbox{for all }t\in\mathbb{R}.
$$
\protect\setcounter{equation}{0}
%
%
%
%
\section{Laurent operators}\label{LO}
This section is the first step of our analysis. 
We will analyse the Laurent operator associated with the symbol $h$ which is much easier 
to study than the Toeplitz operator $T_0$.
\subsection{Notations and main results}
Let  ${\mathcal H}=l^2(\mathbb{Z}, \mathbb{C}^N)$ be the Hilbert
space of square summable vector-valued sequences $ (\gy_n)_{n\in\mathbb{Z}} $ endowed
with its usual scalar product. 
Consider the block Laurent operator associated with the symbol $h$ defined on
   ${\mathcal H}$  by the expression
\begin{equation}
\label{laurent}
(H_0\psi)_n=\sum_{-\infty}^{\infty} A_{n-m}\psi_m, ~~\mbox{
\quad } n\in\mathbb{Z}.
\end{equation}
 In \cite{S3} the case where $h$ is a trigonometric polynomial symbol is studied and applications to different concrete models  are given. Those applications represent additional motivations of our interest to this kind of operators. 
The case where $N=1$ corresponds to the usual scalar Laurent  operators 
that are well studied by different approaches, see 
 \cite{ABC,BS1,BS2} and their references.  

Let us denote by $\textbf{N}$ the multiplication operator defined by 
$(\textbf{N}\psi)_j=j\psi_j$ for all $\psi\in\mathcal{H}$.
Recall that the interpolation space $\mathcal{K}:=(D(\textbf{N}),\mathcal{H})_{1/2,1}$ can be described  by the norm 
$$
\|\psi\|_{1/2,1}=\sum_{j=0}^\infty2^{j/2}\|\theta(2^{-j}|\textbf{N}|)\psi\|
$$ 
where $\theta\in C^\infty_0(\mathbb{R})$ is chosen as in  (\ref{besov}). Clearly,  
$\langle {\textbf{N}}\rangle^{-s}\mathcal{H}\subset\mathcal{K}\subset\mathcal{H}$ for any $s>1/2$. 
 The main result of this section is the following.
\begin{theorem} \label{rbvcontinuity} 
Assume that Assumption \ref{analytique} holds. Then\\
(i) the  spectrum of $H_0$ is purely absolutely continuous outside
$\kappa(h)$ and therefore
$$
\sigma(H_0)=\sigma_{ess}(H_0)=\cup_{i=1}^{i=N}\Sigma_j~,\quad\sigma_p(H_0)\subset\kappa(h)\quad\mbox{ and }~~\sigma_{sc}(H_0)=\emptyset.
$$
(ii)  The  limits 
$
 ({H_0}-x\mp i0)^{-1}:=\lim_{\mu\rightarrow0+}({H_0}-x\mp i\mu)^{-1}
$
exist locally uniformly on  $\mathbb{R}\setminus\kappa(h)$ for the weak* topology of 
$B({\mathcal{K},\mathcal{K}^*})$.\\
(iii) For any $s>1/2$, the maps 
$
 x\mapsto\langle {\textbf{N}}\rangle^{-s}({H}_0-x\mp i0)^{-1}\langle {\textbf{N}}\rangle^{-s}\in B({\mathcal   H})
 $
are locally of class $\Lambda^{s-\frac{1}{2}}$ on  $\mathbb{R}\setminus\kappa(h)$.  \\
(iv) For any  $\varphi\in C^\infty_0(\mathbb{R}\setminus\kappa(h))$ and  $\sigma,\varepsilon>0$.
$$
||\langle {\textbf{N}}\rangle^{-\sigma}e^{-i{H_0}t}\varphi(H_0)\langle {\textbf{N}}\rangle^{-\sigma}||\leq C(1+|t|)^{-\sigma+\varepsilon}.
$$
\end{theorem}
This theorem will be proved in several steps. Let $\mathbb{T}=\mathbb{R}/2\pi\mathbb{Z}$ be the periodized interval $[-\pi,\pi)$.
Consider the discrete Fourier transform
$\mathcal{F}:l^2(\mathbb{Z},\mathbb{C}^N)\rightarrow L^2(\mathbb{T},\mathbb{C}^N)$ defined by
$$
\mathcal{F}(\psi)(p)=\widehat\psi(p)=\sum_{n\in\mathbb{Z}}\psi_ne^{inp}~~,\mbox{ \quad }~~\psi\in l^2(\mathbb{Z},\mathbb{C}^N).
$$
Direct computation shows that,  for each $\psi\in l^2(\mathbb{Z},\mathbb{C}^N)$,
$$
\widehat{H_0\psi}(p)=h(p)\widehat{\psi}(p)~~,\mbox{ \quad }~p\in\mathbb{T}~,
$$
where $h(p)=\sum_{-\infty}^{\infty}A_me^{imp}$ and, therefore, 
 $H_0$ is unitarily equivalent  to the direct integral 
\begin{equation}\label{direct}
\hat{H_0}=\int_{\mathbb{T}}^\oplus h(p)dp\quad\mbox{ acting in }\quad
\hat{\mathcal{H}}=\int_{\mathbb{T}}^\oplus\mathbb{C}^Ndp\;.
\end{equation}
Recall 
 for all $p\in\mathbb{T}$,  $\{\lambda_{j}(p);\,1\leq j\leq N\}$ are the repeated eigenvalues  of
$h(p)$ and $\{W_{j}(p);\,1\leq j\leq N\}$ is a corresponding orthonormal basis
of eigenvectors.
 Hence one may deduce  the following.
\begin{proposition}\label{prop3.1}
The operator $H_0$ is unitarily equivalent to $M=\oplus_{j=1}^{j=N}\lambda_j(p)$ acting in 
the direct sum $\oplus_{j=1}^{j=N}L^2(\mathbb{T})$. In particular, the assertion (i) of Theorem \ref{rbvcontinuity}  follows.
\end{proposition}
\textbf{Proof } Consider the operator 
$
\mathcal{U}:l^2(\mathbb{Z},\mathbb{C}^N)\rightarrow\oplus_{j=1}^{j=N}L^2(\mathbb{T})
$
defined by 
$
\mathcal{U}\psi=(f_j)_{1\leq j\leq N}
$,
where 
\begin{equation}\label{deftrsfF}
f_j(p)=\langle W_j(p),\widehat\psi(p){\rangle}=\sum_{n\in\mathbb{Z}}\langle W_j(p),\psi_n{\rangle}e^{inp}.
\end{equation}
Clearly, $\mathcal{U}$ 
is a unitary operator with the following inversion formula:
$$
(\mathcal{U}^{-1}f)_n=\sum_{j=1}^N\frac{1}{2\pi}\int_{-\pi}^\pi f_j(p)W_j(p)e^{-inp}dp,
$$
for all $f=(f_1,\cdots,f_N)$ in $\oplus_{j=1}^{j=N}L^2(\mathbb{T})$.
Moreover, direct computation shows that,  for each $\psi\in l^2(\mathbb{Z},\mathbb{C}^N)$,
$$
(\mathcal{U}H_0\psi)_j(p)=\langle W_j(p),\widehat{H_0\psi}(p){\rangle}=\langle W_j(p),h(p)\widehat{\psi}(p){\rangle}=\lambda_j(p)f_j(p).
$$
The proof is complete.
%

\subsection{Mourre estimate for $\lambda _{j}$}\label{multiplication}

Here we denote  the operator of multiplication by a function 
$f$ in $L^2(\mathbb{T})$  by the same symbol  $f$ or 
by $f(p)$ when we want to stress the $p$-dependence. 

As $H_0$ is unitarily equivalent to $\oplus_{j=1}^N\lambda_j$,
it is convenient to prove first the  Mourre estimate for each multiplication operator $\lambda_j$.
So, let us  fix $j\in\{1,\cdots,N\}$ and put $\lambda=\lambda_j$. 
Recall that the set  $\kappa (\lambda)$ of critical values of $\lambda$ is finite. 
Let $\Delta$ be a real compact interval, such that ${\Delta}\cap \kappa (\lambda)=\emptyset.$
Clearly, there exists an open interval $I$ containing $\Delta$ and $\overline{I}\cap\kappa (\lambda)=\emptyset$.
In particular,   $\lambda$ is of class $C^\infty$ on $\lambda^{-1}(\overline{I})$
and there exists a constant $c>0$ such that $|\lambda^{\prime }(p)|^2\geq c>0$ on $\lambda^{-1}(\overline{I})$.
Choose a smooth function $\zeta$ supported on $\overline{I}$ and identically equals to 1 on $\Delta$.
Define the function
$$
F(p)=\zeta(\lambda(p))\frac{\lambda'(p)}{|\lambda'(p)|^2}~,\quad p\in\mathbb{T}.
$$
Thus $F\in C^{\infty } ({\mathbb{T}})$, since the support of $\zeta$ does not contain any critical point of $\lambda$.
Moreover
\begin{equation}\label{FF}
F(p) \lambda^{\prime }(p)=\zeta(\lambda(p))~,\quad\forall p\in\mathbb{T}.
\end{equation}  
We are now able to define our conjugate operator $\mathbf{a}$ to $\lambda$ on $L^2(\mathbb{T})$ by
\begin{equation}
\mathbf{a}=  \frac{i}{2}\{F(p) \frac{d}{dp}+ \frac{d}{dp} F(p)\}=
i F(p) \frac{d}{dp} +\frac{i}{2}F^{\prime }(p). \label{ocdef}
\end{equation}
It is clear that $\mathbf{a}$ is an essentially self-adjoint operator  in $L^2(\mathbb{T})$, see for example Proposition 7.6.3 of \cite{ABG}. 
Moreover, direct computation yields to
$$
[\lambda,i\mathbf{a}]=F(p) \lambda^{\prime }(p)=\zeta(\lambda(p)).
$$
So $\lambda$ is of class $C^{1}(\mathbf{a})$,  and by repeating the same calculation, we show that  
$\lambda$ is of class $C^{\infty}(\mathbf{a})$. Moreover,  since the spectral projector of $\lambda$ on $\Delta$ is nothing but the multiplication operator by 
the characteristic function of $\lambda^{-1}(\Delta)$ we have
\begin{eqnarray}\label{mel}
E_{\lambda}(\Delta)[\lambda,i\mathbf{a}]E_{\lambda}{(\Delta)}= E_{\lambda}(\Delta).
\end{eqnarray}
Thus we have proved
\begin{proposition}\label{prop3.3}
The operator $\lambda$ is of class $C^{\infty}(\mathbf{a})$ and $\mathbf{a}$ is  strictly conjugate
to $\lambda$ on $\Delta$.
\end{proposition}
%
%
%
\subsection{Mourre estimate for $H_0$}\label{mestimateH0}
 %
 We start with the following consequence of Proposition \ref{prop3.3}.
\begin{proposition}\label{moureM} 
Let $\Delta$ be a real compact interval such that $\Delta\cap\kappa(h)=\emptyset$. Then
there exists a self-adjoint operator $\widetilde{A}=\widetilde{A}_\Delta$ such that $M$ is of class $C^{\infty}(\widetilde{A})$ and $\widetilde{A}$ is
strictly conjugate to $M$ on $\Delta$. More precisely, we have
$$
E_M(\Delta)[M,i\widetilde{A}]E_M(\Delta)=E_M(\Delta).
$$
\end{proposition}
\textbf{Proof  }  
Recall that  
$
\kappa (h)= \cup_{1\leq j\leq N}\kappa (\lambda _{j})
$ 
is finite. Since $\Delta$ is a real compact interval such that $\Delta\cap\kappa(h)=\emptyset$,
 there exists an open interval $I$ containing $\Delta$ and $\overline{I}\cap\kappa (h)=\emptyset$.
In particular, for any  $j=1,\cdots,N$,  $\lambda_j$ is of class $C^\infty$ on $\lambda_j^{-1}(\overline{I})$
and there exists a constant $c>0$ such that $|\lambda_j^{\prime }(p)|^2\geq c>0$ on $\lambda_j^{-1}(\overline{I})$.
Thus for any $j=1,\cdots,N$,  one may construct the operator $\mathbf{a}_{j}$  using (\ref{ocdef}),  with 
$F=F_{j}$ and $\lambda=\lambda_j$,  that is
\begin{equation}
\mathbf{a}_j= \frac{i}{2}\{F_j(p) \frac{d}{dp}+ \frac{d}{dp} F_j(p)\}=
i F_j(p)\frac{d}{dp} +\frac{i}{2}F_j^{\prime }(p). \label{ocdef1}
\end{equation}
For any $j=1,\cdots,N$,  we have
$$
E_{\lambda_j}(\Delta)[\lambda_j,i\mathbf{a}_j]E_{\lambda_j}{(\Delta)}= E_{\lambda_j}(\Delta).
$$
Now define the self-adjoint operator in 
$%
\oplus _{1\leq j\leq N} L^{2}(\mathbb{T})
$ 
by
\begin{equation}
\widetilde{A}=\oplus _{1\leq j\leq N}\mathbf{a}_{j} .  \label{defdeatld}
\end{equation}
Since, for any Borel set $J$, one has  $E_M(J)=\oplus _{1\leq j\leq N}E_{\lambda_j}(J)$ and 
$[M,i\tilde{A}]=\oplus_{1\leq j\leq N}[\lambda_j,i\mathbf{a}_j]$,
 we immediately  conclude that,
 $$
E_{M}(\Delta)[M,i\tilde{A}]E_{M}{(\Delta)}=\oplus _{1\leq j\leq N}E_{\lambda_j}(\Delta)=E_{M}(\Delta).
$$ 
The proof is complete.

According to the Proposition \ref{prop3.1},  $H_0=\mathcal{U}^{-1}M\mathcal{U}$. Define the operator $\mathbb{A}$ by 
\begin{equation}\label{oct}
\mathbb{A}=\mathcal{U}^{-1}\widetilde{A}\mathcal{U}.
\end{equation}
By a direct application of the previous proposition, we get
\begin{corollary}
Let $\Delta$ be a real compact interval such that $\Delta\cap\kappa(h)=\emptyset$. Then
there exists a self-adjoint operator $\mathbb{A}=\mathbb{A}_\Delta$ such that $H_{0}$ is of class $C^{\infty }(\mathbb{A})$ and $\mathbb{A}$ is
strictly conjugate to $H_{0}$ on $\Delta$:
$$
E_0(\Delta)[H_0,i\mathbb{A}]E_0(\Delta)=E_0(\Delta).
$$
\end{corollary}
In particular, the spectrum of $H_0$ is purely absolutely continuous on $\mathbb{R}\setminus\kappa(h)$.
The last ingredient in the proof of Theorem \ref{rbvcontinuity} is the following Lemma that comes from \cite{S3}.
\begin{lemma}
\label{regularitŽé} 
Let $\mathbb{A}=\mathbb{A}_\Delta$ be a self-adjoint operator defined by (\ref{oct}).
For every integer $m>0$ the operators $\langle {\textbf{N}}\rangle^{-m}\mathbb{A}^{m}$ and
$\mathbb{A}^{m}\langle {\textbf{N}}\rangle^{-m} $ are bounded in $\mathcal{H}$.
\end{lemma}
\textbf{Proof of  Theorem \ref{rbvcontinuity}. }
By combining the last corollary, Lemma \ref{regularitŽé}, and the results of Section  \ref{rappel}, we immediately get our Theorem  \ref{rbvcontinuity}.
\protect\setcounter{equation}{0}
\section{Compact perturbation of $H_0$}\label{perturbŽ}
The second step of our analysis  is to study the preservation of this spectral structure under compact perturbations. 
More specifically, let $V=(V_{ij})_{i,j\in\mathbb{Z}}$ be  a symmetric compact  operator in $\mathcal{H}$.
Hence, by Weyl theorem, $H=H_0+V$ and $H_0$ have the same essential spectra,
$$
\sigma_{ess}(H)=\sigma_{ess}(H_0)=\cup_{j=1}^{j=N}\Sigma_j.
$$
We would like to describe simple conditions on $V$ ensuring the preservation of Theorem \ref{rbvcontinuity} for $H$.
For this let us fix  a real compact interval $\Delta$ such that $\Delta\cap\kappa(h)=\emptyset$ and let
$\mathbb{A}=\mathbb{A}_\Delta$ be  the self-adjoint operator defined by (\ref{oct}).
Our analysis in this section is based on the following two abstract results. 
\begin{theorem}\label{MEFORH} 
Assume that the compact operator $V$ is of class ${C}^{1}(\mathbb{A})$ such that $[V,i\mathbb{A}]$ is compact too. 
Then  there exists  a compact operator $K$ such that
$$
E_H(\Delta)[H,i\mathbb{A}]E_H(\Delta)=E_H(\Delta)+K.
$$
In particular, $\Delta$ contains at most a finite number of eigenvalues of $H$ counted with multiplicities.
 \end{theorem}
 \textbf{Proof. }
 According to the last section, the operator
$H_0$ is of class $C^\infty(\mathbb{A})$. Then, since $H_0$ and $V$ are bounded 
and $V$ is of class $C^{1}(\mathbb{A})$,  one deduces that $H=H_0+V$ is of class  $C^{1}(\mathbb{A})$.  Moreover, 
 for any Borel set $J$,  $E_H(J)-E_{H_0}(J)$ is a compact operator.
Hence, 
\begin{eqnarray*}
E_H(\Delta)[H,i\mathbb{A}]E_H(\Delta)&=&E_H(\Delta)[H_0,i\mathbb{A}]E_H(\Delta)+E_H(\Delta)[V,i\mathbb{A}]E_H(\Delta)\\
&=&E_{H_0}(\Delta)[H_0,i\mathbb{A}]E_{H_0}(\Delta)+K_1\\
&=&E_{H_0}(\Delta)+K_2\\
&=&E_{H}(\Delta)+K_3
\end{eqnarray*}
for some  compact operators $K_1,K_2$ and $K_3$. 
The proof is complete.
 \begin{theorem}\label{rbvcontinuity1}
Assume that $V$ is a self-adjoint compact operator on $\mathcal{H}$ which is of class $\mathcal{C}^{1,1}(\mathbb{A})$.Then
$V$ of class ${C}_u^{1}(\mathbb{A})$ and $[V,i\mathbb{A}]$ is compact so that Theorem \ref{MEFORH}  holds. Moreover,
\begin{enumerate}
 \item $H_0+V$ has no singular continuous spectrum in $\Delta$;
 \item the limits 
 $
 ({H}-x\mp i0)^{-1}:=\lim_{\mu\rightarrow0+}({H}-x\mp i\mu)^{-1}
$
exist locally uniformly on  $\Delta\setminus[\kappa(h)\cup\sigma_p(H)]$ for the weak* topology of 
$B({\mathcal{K},\mathcal{K}^*})$.
\item  In addition, if $V$ is of class $\mathcal{C}^{s+\frac{1}{2}}(\mathbb{A})$, for some $s>1/2$, then  the maps
$$
 x\longmapsto  \langle {\textbf{N}}\rangle^{-s}({H}-x\mp i0)^{-1}\langle {\textbf{N}}\rangle^{-s}\in B({\mathcal   H})
$$
are locally of class  $\Lambda^{s- \frac{1}{2}}$ on $\Delta\setminus[\kappa(h)\cup\sigma_p(H)]$. 
\end{enumerate}
\end{theorem}
\textbf{Proof. } A direct combination of Theorem \ref{MEFORH} , Lemma \ref{regularitŽé}, and the results of Section  \ref{rappel}.

Now the efficiency of these theorems in applications lies on our ability to verify their regularity hypothesis for concrete $V$'s.
We start with (see \cite{ABG} for a continuous version) the following  auxiliary result.
\begin{lemma}
Let $V$ be a bounded operator on $\mathcal{H}$ and $\xi\in C^\infty(\mathbb{R})$ such that $\xi(x)=0$  around 0 and $\xi(x)=1$ for $|x|>a$, for some $a>0$. Then the  operator $V$ is compact on $\mathcal{H}$ if and only if
$\lim_{r\rightarrow+\infty}\|\xi(\mathbf{N}/r) V\|=0$.
\end{lemma} 
\textbf{Proof }  First, it is obvious that 
$$
\slim_{r\rightarrow+\infty}\xi(\mathbf{N}/r)=0.
$$
Hence, if  $V$ is compact then $\lim_{r\rightarrow+\infty}\|\xi(\mathbf{N}/r) V\|=0$.

Reciprocally, since $\zeta=1-\xi\in C^\infty(\mathbb{R})$ with a compact support, $V_k=\zeta(\mathbf{N}/k) V$ is a compact operator, for any integer $k$.
Since, $\|V-V_k\|=\|\xi(\mathbf{N}/k) V\|$ tends to zero at infinity we get the desired result.

\begin{example}
Let $f$ be a the  self-adjoint  matrix-valued continuous function defined on $\mathbb{T}$   and let 
$\{\hat{f}_n\}_{n\in\mathbb{Z}}$ be its Fourier coefficients. Consider the operator 
$W=W(f)$ defined by
$$
W_{ij}=
\left\{
\begin{array}{ccc}
\hat{f}_{i-j}&  \mbox{if} &i\geq0,j\leq-1,   \\
 \hat{f}_{i-j}&    \mbox{if} &i\leq-1,j\geq0, \\
   0& \mbox{if} & \mbox{ otherwise. }
\end{array}
\right.
$$
Then $W(f)$ is a compact operator. Indeed, if $f$ is a trigonometric polynomial then $W(f)$ is of finite rank, is clear from the fact that 
$$
\|\xi(\langle \mathbf{N}\rangle/r) W\|=0 ~,\quad\mbox{for $r$ sufficiently large}.
$$
Now, the proof can be finished by the celebrated Fejer's Theorem. 
This example is related to Hartman's theorem on the compactness of  Hankel operators \cite{H}. 
\end{example} 
\begin{corollary}
Let $V$ be a bounded symmetric operator in $\mathcal{H}$ such that
$$
\lim_{r\to+\infty}r\|\xi(\mathbf{N}/r)V\|=0.
$$
Then  outside  $\kappa(h)$ the eigenvalues of $H=H_0+V$ are all finitely degenerate and 
 their possible accumulation points are included in $\kappa(h)$.
\end{corollary}
\textbf{Proof. } It suffices to show that for any real compact interval  $\Delta$  such that $\Delta\cap\kappa(h)=\emptyset$, 
 Theorem \ref{MEFORH} applies, so that $\Delta$  contains
  at most a finite number of eigenvalues of $H$ counted with multiplicities.
Under our condition the operators $V$ and $\langle\textbf{N}\rangle V$ are compact on $\mathcal{H}$. 
Let $\mathbb{A}=\mathbb{A}_\Delta$ be  the self-adjoint operator defined by (\ref{oct}).
According to Lemma \ref{regularitŽé} we get 
$\mathbb{A} V=\mathbb{A}\langle\textbf{N}\rangle^{-1}\cdot\langle\textbf{N}\rangle V$ is compact.
Hence $H$  is of class $C^{1}(\mathbb{A})$ such that $[V,i\mathbb{A}]$ is compact. 
The proof is complete.

\textbf{Remark }
According to Theorem 6.1 of \cite{BS1}, $V$ is of class $\mathcal{C}^{1,1}(\mathbb{A})$ if
\begin{equation}\label{C1,1}
\int_1^{+\infty}\|\theta(\langle \mathbf{N}\rangle/r) V\|dr<\infty.
\end{equation}
Recall that $\theta\in C^\infty_0(\mathbb{R})$ with  $\theta(x)>0$ if $1<x<2$ and $\theta(x)=0$ otherwise.
Moreover,  if
\begin{equation}\label{Cs}
\sup_{r\geq 1}\|r^s\theta(\langle \mathbf{N}\rangle/r) V\|<\infty
\end{equation}
then $V$ is of class $\mathcal{C}^{s}(\mathbb{A})$.

Hence to show that $V$ is of some regularity class with respect to $\mathbb{A}$  it is enough to have a good estimation of the norms  $\|\theta(\langle \mathbf{N}\rangle/r) V\|$.
For we need the following notations. For an operator $B=(B_{ij})_{i,j\in\mathbb{Z}}$ defined formally on $\mathcal{H}$ by
$
(B\psi)_i=\sum_{j=-\infty}^\infty B_{ij}\psi_j
$.
It is known and easy to show that $B$ is a bounded operator in $\mathcal{H}$ if  
\begin{eqnarray}
r(B)&:=&\sup_{i\in\mathbb{Z}}\sum_{j\in\mathbb{Z}}\|B_{ij}\|<\infty,\\
c(B)&:=&\sup_{j\in\mathbb{Z}}\sum_{i\in\mathbb{Z}}\|B_{ij}\|<\infty.
\end{eqnarray}
In this case,
$$
\|B\|\leq\sqrt{r(B)c(B)}.
$$
This is the so-called Schur test. Notice that if $B$ is symmetric then  $r(B)=c(B)$ and $\|B\|\leq c(B)$.
The proofs are similar to the scalar case and will be omitted.

We will also use the following Hilbert-Schmidt test. Recall that $B$ is Hilbert-Schmidt if
$$
\|B\|_2=\sqrt{\sum_{i,j\in\mathbb{Z}}\|B_{ij}\|^2}<\infty.
$$
In this case $B$ is bounded and $\|B\|\leq \|B\|_2$.
Here also the proof is similar to the scalar case.
\begin{example} Consider  the operators $B$ and $B'$ on $\mathcal{H}$ such that the matrix entries:
$$
\|B_{ij}\|=\frac{1}{(1+|i-j|)^2}\quad \mbox{and}\quad \|B'_{ij}\|=\frac{1}{(1+|i|)(1+|j|)}~.
$$
Then clearly Schur test shows that $B$ is bounded while the Hilbert-Schmidt test  fails.
In contrast  Schur  test breaks down for $B'$, while Hilbert-Schmidt test applies. 
\end{example}
\begin{proposition}\label{perturbation-regularité1}
Let $V=(V_{ij})_{i,j\in\mathbb{Z}}$ be a bounded symmetric operator in $\mathcal{H}$. Put 
\begin{eqnarray}
n_V(r)&:=&(\sup_{r\leq|i|\leq 2r }\sum_{j\in\mathbb{Z}}\|V_{ij}\|)\cdot(\sup_{ j\in\mathbb{Z}}\sum_{r\leq|i|\leq 2r}\|V_{ij}\|),\\
p_V(r)&:=&\sum_{\substack{r\leq|i|\leq 2r \\ j\in\mathbb{Z}}}\|V_{ij}\|^2.
\end{eqnarray}
 Then $V$ is of class $\mathcal{C}^{1,1}(\mathbb{A})$ if
$
\int_1^{\infty}\sqrt{n_V(r)}dr<\infty
$,
respectively, if
$
\int_1^{\infty}\sqrt{p_V(r)}dr<\infty
$.
\end{proposition}
\textbf{Proof. } According to the preceding discussion if we use the Schur (respectively Hilbert-Schmidt) test one may estimate
$
\|\theta(\langle \mathbf{N}\rangle/r) V\|
$
by $\sqrt{n_V(r)}$ (respectively $\sqrt{p_V(r)}$).
Hence the hypothesis of the proposition implies 
the estimate (\ref{C1,1}).  The proof is complete.
\begin{corollary}
Assume that $V=V_{1}+V_{2}$ so that 
$$
\int_1^{+\infty}\sqrt{n_{V_1}(r)}+\sqrt{p_{V_{2}}(r)}dr<\infty.
$$
Then $V$ is of class $\mathcal{C}^{1,1}(\mathbb{A})$
\end{corollary}
Similarly we have
\begin{corollary}\label{perturbation-regularité}
Let $V=V_{1}+V_{2}$ be a bounded symmetric operator in $\mathcal{H}$ such that 
$$
\sup_{r\geq1}r^s(\sqrt{n_{V_1}(r)}+\sqrt{p_{V_{2}}(r)})<\infty.
$$
Then $V$ is of class $\mathcal{C}^s(\mathbb{A})$.
\end{corollary}
\textbf{Remark } Notice also that if 
$\langle\textbf{N}\rangle^{s}V\in B(\mathcal{H})$, for some
 $s>0$ then $V$ is of class $\mathcal{C}^s(\mathbb{A})$.  In particular, this holds if $\|\langle\textbf{N}\rangle^{s}V\|_2<\infty$
or $
r(\langle\textbf{N}\rangle^{s}V)\cdot c(\langle\textbf{N}\rangle^{s}V)<\infty
$. Indeed, according to the criteria given in the appendix of \cite{BS1} it is enough to show that the product
$<\mathbb{A}>^{s}V$ is bounded in $\mathcal{H}$. But the lemma \ref{regularitŽé}  shows that
$\langle {A}\rangle^{s }\langle {\textbf{N}}\rangle^{-s }$ is bounded. Thus 
 $
<\mathbb{A}>^{s}V=(\langle {A}\rangle^{s }\langle {\textbf{N}}\rangle^{-s })\langle {\textbf{N}}\rangle^{s }V
 $
is  bounded.
\begin{example} Assume that $V$ is compactly supported, that  is $V_{ij}=0$ for all $|i|>M$ or $|j|>M$ for some integer $M$.
Then $V$ is of finite rank and of class $C^\infty(\mathbb{A})$.  Indeed, by the Corollary \ref{perturbation-regularité}, 
we see that $V$ is of class $\mathcal{C}^s(\mathbb{A})$ for any $s>0$.
\end{example}
\begin{example} Let $\psi\in\mathcal{H}$ and consider the rank one operator $V=\langle\cdot,\psi\rangle\psi$.
One may show that $\|\theta(\mathbf{N}/r) V\|=\|\psi\|\cdot\|\theta(\mathbf{N}/r)\psi\|$.
Hence the regularity of $V$ can be described in terms of the regularity of $\psi$. For example,
$V\in\mathcal{C}^{1,1}$ (respectively $V\in\mathcal{C}^{s}, s>0$) if and only if $\psi\in\mathcal{H}_{1,1}$ (respectively
$\psi\in\mathcal{H}_{s,\infty}$). Here $\mathcal{H}_{s,p}$ are the Besov spaces associated to the operator $\mathbf{N}$, see \cite{ABG}.
\end{example}
\begin{example}\label{exponentialdecay}
Assume that there exist a constant $C>0$ and $\kappa>0$ such that
$$
\|V_{ij}\|\leq Ce^{-\kappa(|i|+|j|)}.
$$
Then, according to the Corollary \ref{perturbation-regularité}, $V$ is of class $C^\infty(\mathbb{A})$.
\end{example}
\begin{example}
Assume that there exist a constant $C>0$ and $s>0$ such that
$$
\|V_{ij}\|\leq C(1+|i|+|j|)^{-1-s}.
$$
Then  the Corollary \ref{perturbation-regularité}  and direct calculation show that $V$ is of class $\mathcal{C}^s(\mathbb{A})$. 
\end{example}
\begin{example}
Assume that there exist a constant $C>0$ and $\sigma>1/2$ such that
$$
\|V_{ij}\|\leq C\frac{1}{(1+|i|)^\sigma(1+|j|)^\sigma}.
$$
Then  clearly
\begin{eqnarray*}
p_V(r)&=&\sum_{r\leq|i|\leq 2r }\sum_{ j\in\mathbb{Z}}\|V_{ij}\|^2\\
&\leq&C\sum_{r\leq|i|\leq 2r }\sum_{ j\in\mathbb{Z}}\frac{1}{(1+|i|)^{2\sigma}(1+|j|)^{2\sigma}}
\leq C_\sigma\sum_{r\leq|i|\leq 2r }\frac{1}{(1+|i|)^{2\sigma}}
\end{eqnarray*}
The last term behaves like $r^{-2\sigma+1}$ as $r$ goes to infinity. Hence, $V$ is of class $C^{\sigma-\frac{1}{2}}(\mathbb{A})$. Notice that if $\sigma\leq1$ then one may not use $n_V(r)$ because it diverges.
\end{example}

We close this section by some scattering consideration. Let  $\mathcal{K}^{*\circ}$  be the closure of $\mathcal{H}$ in  
$\mathcal{K}^*$.
In fact it is enough to see that our Theorems \ref{moureM} , \ref{rbvcontinuity1} and \ref{MEFORH} 
allows us to apply the Proposition 7.5.6, page 323 of \cite{ABG} to get the following.

\begin{corollary}\label{waveoperators1}
Assume that  $V$ satisfies (\ref{C1,1}) and extends to a bounded operator from $\mathcal{K}^{*\circ}$ to 
 $\mathcal{K}$.  Let $P_{ac}(H_0)$ be the spectral projector of $H_0$ on its purely absolutely continuous component.
 Then the wave operators
 $$
 \Omega^\pm(H,H_0)=\slim_{t\rightarrow\pm\infty}e^{iHt}e^{-iH_0t}P_{ac}(H_0)
 $$
exist and they are complete, i.e., their range equal $\mathcal{H}_{ac}(H)$ (the absolutely continuous subspace of $H$).
\end{corollary}
\begin{example}  \label{periodique}
Let $H_0$ be the Laurent operator associated with the symbol 
\begin{equation}\label{h(p)periodique}
h(p)=\left(
\begin{array}{ccccc}
  b_1&a_1&0&\cdots&a_Ne^{ip}  \\
  a_1&b_2&a_2&0&\cdots\\
  0&\ddots&\ddots&\ddots&\vdots\\
  0&\cdots&a_{N-2}&b_{N-1}&a_{N-1}\\
  a_Ne^{-ip}&\cdots&0&a_{N-1}&  b_N  
\end{array}
\right) ~,\quad  p\in\mathbb{T} ,
\end{equation}
with $a_i>0, b_i\in \mathbb{R}$. 
Then one may prove that we have, see \cite{S3} for more details,
\begin{enumerate}
\item there exist $\alpha_1<\beta_1\leq\alpha_2<\beta_2\cdots\leq\alpha_N<\beta_N$  such that  the spectrum of $H_0$ is purely absolutely continuous and
$
\sigma(H_0)=\cup_{j=1}^N[\alpha_j,\beta_j].
$
\item 
There exists a conjugate operator $\mathbb{A}$ for $H_0$  i.e.
$
\mu^{\mathbb{A}}(H_0)=\mathbb{R}\setminus\{\alpha_1,\beta_1,\alpha_2,\beta_2\cdots,\alpha_N,\beta_N\}$. 
\item  Theorems \ref{rbvcontinuity1} and its corollaries can be used in this context to get a large class of compact perturbation $V$
that ensures the absence of singular continuous spectrum of $H=H_0+V$ as well as the existence and completeness of the wave operators, but we will not do this separately.
\end{enumerate}
Note that $H_0$ is unitarily equivalent to the scalar periodic  Jacobi operator acting in $l^2(\mathbb{Z},\mathbb{C})$ by
$$
(J_0 \psi)_n=a_n\psi_{n+1}+b_n\psi_n+a_{n-1}\psi_{n-1}, ~~ n\in\mathbb{Z}
$$
such that $a_n>0$ and  $b_n\in\mathbb{R}$ with $a_{j+N}=a_j$ and $b_{j+N}=b_j$. One may also see  \cite{DKS,KvM,vM} and references therein for related results on this model. 
\end{example}
\begin{example} Notice that if in the preceding example one of the $a_i$'s is zero, say $a_1=0$, then
all the spectral bands of $H_0$ are degenerate $\alpha_1=\beta_1, \alpha_2=\beta_2\cdots,\alpha_N=\beta_N$.
Therefore $\sigma_{ess}(H_0)=\{\alpha_1,\alpha_2,\cdots,\alpha_N\}$.
\end{example}

%
%
%
\section{Proof of Theorems \ref{thm1} and \ref{thm2}}\label{toelplitz-analytique}
%
Recall that $\mathbb{Z}_+=\{0,1,\cdots,\}$ and $\mathbb{Z}_-=\{\cdots,-2,-1\}$ are the sets of nonnegative and negative integers. Of course  $\mathbb{Z}=\mathbb{Z}_-\cup\mathbb{Z}_+$ the set of all integers. Let  ${\mathcal H}^\pm=l^2(\mathbb{Z}_\pm, \mathbb{C}^N)$, for some integer  $N\geq1$, be the Hilbert space of square summable vector-valued sequences 
$ (\gy_n)_{n\in\mathbb{Z}_\pm} $.   In the sequel ${\mathcal H}^\pm$ will be  canonically identified with a subspace of ${\mathcal H}$ so that
 ${\mathcal H}={\mathcal H}^-\oplus{\mathcal H}^+$.
With these notations the operator $H_0$ decomposes as
$$
H_0=\left(
\begin{array}{ccc}
 H_0^- & H_0^\mp   \\
 H_0^\pm&  H_0^+
\end{array}
\right).
$$
Of course, the operators $T_0,H_0^+$ and $ H_0^-$ have the same spectral properties  and $H_0^\pm$ are Hankel operators.
Moreover,
$$
H_1=\left(
\begin{array}{ccc}
 H_0^- &0     \\
0&  H_0^+
\end{array}
\right)=H_0-\left(
\begin{array}{ccc}
 0 &  H_0^\mp  \\
 H_0^\pm&   0
\end{array}
\right).
$$
Hence, to study $T_0$ it is more convenient to study the direct sum 
$H_1$,  since it is a perturbation  of $H_0$ by the operator
$$
W=\left(
\begin{array}{ccc}
 0 &  H_0^\mp  \\
 H_0^\pm&   0
\end{array}
\right).
$$
\begin{example}
Assume that $A_j=0$ if $|j|>M$ for some integer $M\geq1$. This means that the symbol $h(p)=\sum_{j=-M}^MA_{j}e^{ijp}$
is a trigonometric polynomial.  Then $W$ is a finite rank operator since 
$$
W_{ij}=
\left\{
\begin{array}{ccc}
  A_{i-j}&  \mbox{if} &i\geq0,j\leq-1,   \\
  A_{i-j}&    \mbox{if} &i\leq-1,j\geq0, \\
   0& \mbox{if} & \mbox{ otherwise }
\end{array}
\right.=
\left\{
\begin{array}{ccc}
  A_{i-j}&  \mbox{if} &i\geq0,i-M\leq j\leq -1,   \\
  A^*_{j-i}&    \mbox{if} &j\leq 0,j-M\leq i\leq-1, \\
   0& \mbox{} & \mbox{ otherwise. }
\end{array}
\right.
$$
\end{example}
\textbf{Proof of the Corollary \ref{corollaire1}}. By a direct combination between the last
 example with our Theorem \ref{MEFORH}, we get all the corollary except the finiteness of the number 
of the eigenvalues counted with multiplicities. But it is a general result that a finite rank perturbation may create at most a finite number of eigenvalues outside the essential spectrum.

\textbf{Proof of Theorems \ref{thm1} and \ref{thm2}}.
Here again it is enough to study the perturbation $W$ in order to apply 
our Theorems \ref{MEFORH} and \ref{rbvcontinuity1} and the preceding discussion.
We have
$$
W_{ij}=
\left\{
\begin{array}{ccc}
  A_{i-j}&  \mbox{if} &i\geq0,j\leq-1,   \\
  A_{i-j}&    \mbox{if} &i\leq-1,j\geq0, \\
   0& \mbox{if} & \mbox{ otherwise. }
\end{array}
\right.
$$
It is enough to show that $W$ is covered by the Example \ref{exponentialdecay}. 
But this is clear since the symbol $h$ is analytic so that its Fourier coefficients $A_j$ decay exponentially. Therefore $W$ is of class $C^\infty(\mathbb{A})$ and the proof is complete.
%
%
%
\section{Proof of Theorems \ref{toeplitz-perturbed01} and \ref{toeplitz-perturbed02}}\label{Toeplitz-pertubation}
%
%
In this part our goal is to  study $T=T_0+v$ acting in $\mathcal{H}^+$. Here again it is more convenient to study the operator
$H$ acting in $\mathcal{H}=\mathcal{H}^-\oplus\mathcal{H}^+$ giving by
$$
\left(
\begin{array}{ccc}
 H_0^- &0     \\
0&  H_0^++v
\end{array}
\right)=H_1+\left(
\begin{array}{ccc}
 0 &  0   \\
0&   v
\end{array}
\right).
$$
The first term in the right hand side is of class $C^\infty(\mathbb{A})$ and $\mathbb{A}$ is locally conjugate to 
$H_1$ on $\mathbb{R}\setminus\kappa(h)$. It is enough to study the operator $V$ obtained from $v$ by
$$
V_{ij}=\left\{\begin{array}{ccc}
v_{ij}&\mbox{ if } i,j\geq0,   \\
0&  \mbox{ otherwise.} 
\end{array}
\right.
$$
But then our conditions (\ref{admissible1}) (respectively (\ref{admissible2})) on $v$ implies directly that the corresponding $V$ 
is of class $\mathcal{C}^{1,1}(\mathbb{A})$ (respectively $\mathcal{C}^{s+\frac{1}{2}}(\mathbb{A})$).  Theorems \ref{MEFORH} and \ref{rbvcontinuity1}
allows us to finish the proof exactly in the same way as in the preceding section.

%

\textbf{Proof of Corollary \ref{waveoperators}}.
A direct application of our Theorems \ref{thm1}, \ref{thm2} and \ref{toeplitz-perturbed01} 
allows us to apply Proposition 7.1.5, page 282 of \cite{ABG}. Indeed, the space $\mathcal{K}_+$ is of cotype 2 and has 
the bounded approximation property, see Lemma 2.8.6 and Theorem 2.8.7, page 68 of \cite{ABG}.
\begin{example} 
Let $T_0$ be the Toeplitz operator associated with the symbol $h(p)$ given by
(\ref{h(p)periodique}). 
Then we have 
$$
\sigma_{ac}(T_0)=\sigma_{ess}(T_0)=\sigma_{ess}(H_0)=\cup_{j=1}^N[\alpha_j,\beta_j].
$$
Moreover, $T_0$ has at most finite number of eigenvalues in the spectral gaps. In fact, according to \cite{C}, $T_0$ has at most one simple eigenvalue in each gap.
In \cite{C}  (see also \cite{Co})  sufficient conditions on a perturbation $V$ that ensure the finiteness of eigenvalues of $T=T_0+V$  in the spectral gaps are given.
Here Theorems \ref{toeplitz-perturbed01} and \ref{toeplitz-perturbed02} complete the picture. Indeed they can be used in this context to get a large class of compact perturbation $v$
that ensures the absence of singular continuous spectrum of $T=T_0+v$ as well as the existence and completeness of the wave operators, but we will not do this separately.
Finally, note that $T_0$ is unitarily equivalent to the scalar periodic  Jacobi operator defined in $l^2(\mathbb{Z}_+,\mathbb{C})$ by
$$
(J_0^+ \psi)_n=a_n\psi_{n+1}+b_n\psi_n+a_{n-1}\psi_{n-1}, ~~ n\geq0
$$
such that $\psi_{-1}=0, a_n>0$ and  $b_n\in\mathbb{R}$ with $a_{j+N}=a_j$ and $b_{j+N}=b_j$. 
\end{example}
%
%

%

\begin{thebibliography}{99}
%
\bibitem{ABG}
{W.  Amrein , A. Boutet de Monvel  and V. Georgescu,
} { C$_0$-Groups, Commutator Methods and Spectral
Theory of $N$-Body Hamiltonians\/},
Birkh{\"a}user, Progress in Math.\ Ser., 135, Basel.
%
\bibitem{ABC} M.A. Astaburuaga, O. Bourget, V.H. Cort\'es, Spectral stability for compact perturbations of Toeplitz matrices,
 J. Math. Anal. Appl. 440 (2016), no. 2, 885-910. 
%
\bibitem{Ba}  H. Baumg{\"a}rtel, Analytic perturbation theory for matrices and operators. Operator Theory: Advances and Applications, 15. Birkh{\"a}user Verlag, Basel, 1985.
%
\bibitem{B}
 Yu. M. Berezanskii,
{ Expansions in eigenfunctions of selfadjoint operators\/}, Vol. 17,
Translation of Mathematical Monographs,
AMS, Providence, RI, 1968.
%
\bibitem{BoS}{T. Bouhennache, J.  Sahbani,} Conjugate operators from dispersion curves for perturbations of fibered systems in cylinders. Comm. Partial Differential Equations 26 (2001), no. 7-8, 1091-1115.
%
\bibitem{BoSi} A. B\H{o}ttcher, B.  Silbermann,  Analysis of Toeplitz operators. Second edition. Prepared jointly with Alexei Karlovich. Springer Monographs in Mathematics. Springer-Verlag, Berlin, 2006. 
\bibitem{BG}
{ A. Boutet de Monvel, V. Georgescu,}
Boundary values of the resolvent of a self-adjoint operator: higher
order estimates.  Algebraic and geometric methods in mathematical
physics,  9-52, Math. Phys. Stud., 19, Kluwer
Dordrecht, 1996.
%
 %
\bibitem{BGS2}
{ A. Boutet de Monvel, V. Georgescu and J. Sahbani, }
Boundary values of regular resolvent families.  Helv. Phys. Acta
71  (1998), 518--553.
%
\bibitem{BS1} {A. Boutet de Monvel, J. Sahbani,},
  On the spectral properties of discrete Schr{\"o}dinger operators:
The multidimensional case. Rev. Math. Phys., 11 (1999),
1061--1078.
%
\bibitem{BS2} {A. Boutet de Monvel, J. Sahbani,}
 Anisotropic Jacobi matrices with absolutely continuous spectrum. C. R. A. S. Paris SŽr. I Math. 328 (1999), no. 5, 443-448.
 %
  \bibitem{CR}  K. Clancey, L.  Rodman,  Eigen-spaces of families of Toeplitz operators with rational matrix symbols. Manuscripta Math. 45 (1983), no. 1, 1-12.
 \bibitem{CH}  S. M. Cheung, H. Hochstadt,  An inverse spectral problem. Linear Algebra and Appl. 12 (1975), no. 3, 215-222.
%
 \bibitem{Co} P. A. Cojuhari,  Spectrum of the perturbed matrix Wiener-Hopf operators. Mat. Issled., no.61 (1981), 29-39 (Russian); MR 616842 (82 j: 47040).
%
 \bibitem{C1} P. A. Cojuhari, On the spectrum of a class of block Jacobi matrices. Operator theory, structured matrices, and dilations, 137-152, Theta Ser. Adv. Math., 7, Theta, Bucharest, 2007.
 %
  \bibitem{C} P. A. Cojuhari,  Discrete spectrum in the gaps for perturbations of periodic Jacobi matrices. J. Comput. Appl. Math. 225 (2009), no. 2, 374-386.
%
  \bibitem{CF} P. A. Cojuhari, I. Feldman,  
  On the finiteness of the discrete spectrum of a perturbed 
  self-adjoint operator, Mat. Issled, (47)(1978), 35-40 (Russian). 
  %
   \bibitem{DKS} D.  Damanik, R. Killip and  B. Simon,  Perturbations of orthogonal polynomials with periodic recursion coefficients. Ann.  Math. (2) 171 (2010) no. 3, 1931-2010.
  %
  \bibitem{D} R.G. Douglas, Banach Algebra Techniques in Operator Theory, Second edition. 
  Graduate Texts in Mathematics, 179. Springer-Verlag, New York, 1998.
  %
 \bibitem{F} I. Feldman, Finiteness of the discrete spectrum of some block Toeplitz operators, Integral Equations Operator Theory \textbf{16}(1993), 385-391.
 %
 %
\bibitem{GK} I. C. Gohberg , M. G. Krein,  Systems of integral equations on a Half line with kernels depending on the difference of arguments, Amer. Math. Soc. Transl. Ser.2 
\textbf{14}(1960), 217-287.
%
\bibitem{GF} I. C. Gohberg , I. A. Feldman,  Convolution Equations and Projection Methods for Their Solution, Transl. Math. Monographs, vol.41, Amer. Math. Soc., Providence, RI 1974.
%
\bibitem{H} Ph. Hartman, On completely continuous Hankel matrices, Proc. Amer. Math. Soc. 9 (1958) 862-866.
%
\bibitem{HW} Ph. Hartman,  A. Wintner, The spectra of Toeplitz's matrices, Amer. J. Math. 76 (1954) 867-882.
%
 \bibitem{H} G. Heinig, On invertibility and spectra of matrix Wiener-Hopf operators, Math. Sb. 91\textbf{133}(1973), 253-266.
  %
 \bibitem{Ka} T. Kato, Perturbation theory for linear operators. Springer-Verlag, 1966.
%
\bibitem{KvM}  M. Kac, P. van Moerbeke,  On some periodic Toda lattices. Proc. Nat. Acad. Sci. U.S.A. 72 (1975), 1627-1629.
%
%
\bibitem{K}  {M. G. Krein, }  Integral equations on a half-line with kernel depending upon the difference of the arguments, Uspehi Math. Nauk, 13:5 (1958), 3-120 (Russian); Amer. Soc. Transl. 22 (2) (1962), 163-288.
%
\bibitem{vM} P. van Moerbeke, The spectrum of Jacobi matrices. Invent. Math. 37 (1976), no. 1, 45-81. 
%
\bibitem{M} E. Mourre, Absence of singular continuous spectrum for certain self-adjoint operators. Comm. Math. Phys. 78 (1980/81), no. 3, 391-408.
%
%
\bibitem{Pe}  V. V. Peller,  Hankel operators and their applications. Springer Monographs in Mathematics. Springer-Verlag, New York, 2003. 
\bibitem{Pu} C.R. Putnam, On Toeplitz matrices, absolute continuity, and unitary equivalence.
Pacific J. Math. 9 (1959) 837-846. 
%
 \bibitem{R} L. Rodman,  On the structure of self-adjoint Toeplitz operators with rational matrix symbols. Proc. Amer. Math. Soc. 92 (1984), no. 4, 487-494.
%
 \bibitem{Ros} M. Rosenblum, The absolute continuity of Toeplitz's matrices.
Pacific J. Math. 10 (1960) 987-996. 
%
\bibitem{S1}  J. Sahbani,
{ The Conjugate Operator Method For Locally Regular Hamiltonians,}
J. Operator Theory, 38 (1997), 297--322.
%
%
%
\bibitem{S3}
{ J. Sahbani,}  Spectral theory of certain block Jacobi matrices and applications, J. Math. Anal. Appl.  J. Math. Anal. Appl. 438 (2016), no. 1, 93-118. 
  
%
%
\bibitem{Z} A. Zygmund, Trigonometric Series, Cambridge University Press.
\end{thebibliography}
\end{document}